\input amstex
\documentstyle{amsppt}
\NoRunningHeads
\TagsOnRight
\def\N{\Bbb N}
\def\RR{\Bbb R}

\def\NC{N_\Bbb C}
\def\prlim{\mathop{\roman{pr\ lim}}}
\def\indlim{\mathop{\roman{ind\ lim}}}
\def\H{\Cal H}
\def\l{\langle}
\def\r{\rangle}
\def\L{\langle\!\langle}
\def\R{\rangle\!\rangle}

\document

We study the gamma noise measure $\mu$ on Borel sets of Schwartz space
$S'(\RR)$, which turns out to be an interesting example of compound
Poisson measure. The orthogonal analysis with respect to $\mu$ relates to
biorthogonal one analogously to the case of Poisson measure. This fact
enables to solve stochastic Wick--Skorokhod equations involving gamma
noise.

1. Let $N'\supset\H\supset N$ be a nuclear rigging of the real Hilbert
space $\H=\roman L_{2,Re}(\RR)$ by the real Schwartz space
$N:=S(\RR)=\prlim_{p\in\N}\H_p$ and
$N':=S'(\RR)=\indlim_{p\in\N}\H_{-p}$. Here the family
$\{\H_p,\ p\in\N\}$ of the Hilbert spaces satisfies the condition:
$\forall p\in\N$ $\exists p'\in\N:\ \H_{p'}\hookrightarrow\H_p$ and the
embedding belongs to Hilbert--Schmidt class; $\H_{-p}$ is the negative
space of the chain $\H_{-p}\supset\H\supset\H_p$ (see, e.g., [2]). We
denote by $|\cdot|$ and $|\cdot|_p$ the norms in $\H$ and $\H_p$,
$p\in\Bbb Z\backslash\{0\}$, respectively and by $\l\cdot,\cdot\r$ the dual
pairings between $\H_p$ and $\H_{-p}$ and between $N$ and $N'$, given by
the extension of the inner product in $\H$ and preserve these notations for
complexifications and tensor powers of abovementioned spaces.

We define the gamma noise measure $\mu$ on $(N',\Cal B(N'))$ by its
characteristic functional
$$
C_\mu(\theta)=\int_{N'}\exp\{i\l x,\theta\r\}\mu(\roman dx)
=\exp\{-\l 1,\log(1-i\theta)\r\},\ \theta\in N \tag1
$$
via the Bochner--Minlos theorem. This name is natural for $\mu$, because
$C_\mu$ coincides with the characteristic functional of the gamma noise,
i.e. generalized stochastic process with independent values
$\xi'=\{\xi_t',\ t\in\RR\}$ being the distributional derivative of gamma
process $\xi$ (for more details on generalized stochastic processes see
[5], on gamma process---[13]).

Indeed, setting in \thetag1 $\theta=\lambda\chi_{[0\wedge t,0\vee t]}$,
$\lambda\in\RR$, $t\in\RR\backslash\{0\}$ one obtain
$C_\mu(\theta)=(1-i\lambda)^{-|t|}$. This function coincides with Fourier
transform of the density $p_{|t|}(s)=s^{|t|-1}e^{-s}\chi_{(0,\infty)}(s)/
\Gamma(|t|)$, $s\in\RR$ of $\Gamma(|t|)$--distribution. Consider the
probability space $(N',\Cal B(N'),\mu)$. The random variable $\xi_t(x)
:=\l x,\chi_{[0\wedge t,0\vee t]}\r$ has $\Gamma(|t|)$ distribution, so the
family of random variables $\xi=\{\xi_t,\ t\in\RR\}$ ($\xi_0:=0$) is
gamma process.

The gamma process is a separable stochastic process with independent
increments, so the characteristic function
$\bold E\exp\{i\lambda(\xi_{t_2}-\xi_{t_1})\}$, $t_1<t_2$ of its increment
has L\'evy--Khinchine representation (see, e.g., [6,10]). By \thetag1 this
representation has the following form
$$
\int_{N'}\exp\{i\lambda\l x,\chi_{[t_1,t_2]}\r\}\mu(\roman dx)
=\exp\big\{(t_2-t_1)\int_0^\infty(e^{iu\lambda}-1)\frac{e^{-u}}u
\roman du\big\}. \tag2
$$
In other words, the equality \thetag2 means that the L\'evy triple
$(a,\sigma^2,\beta)$ of the gamma process $\xi$ coincides with
$$
a=\int_0^\infty\frac{e^{-u}}{1+u^2}\roman du,
\ \ \sigma^2=0,
\ \ \beta(\Delta)=\int_{\Delta\cap(0,\infty)}\frac{e^{-u}}u\roman du.
$$
Thus, the gamma process is a non-decreasing compound Poisson process such
that almost every sample path has infinitely many jump discontinuities in
each open interval and is locally bounded (see, e.g., [6]). This fact,
equality \thetag1 and some results from [13, ch. 3] enables us to describe
some sets of full $\mu$-measure. Below we denote by $\frak D(\RR)$ the
collection of all the locally countable sets $\frak d\subset\RR$.

\proclaim{Theorem 1}
1) There exists $p\in\RR$ such that $\mu(\H_{-p})=1$.

2) The set
$$
M:=\big\{\sum_{t_k\in\frak d}h_k\delta_{t_k}:\ h_k>0,
\ \frak d\in\frak D(\RR)\big\}\cap C_0(\RR)'\subset N'
$$
has the full $\mu$-measure.

3) Denote for $y\in M$ $y(\tau)=\l y,\chi_{[0\wedge\tau,0\vee\tau]}\r$.
Then
\roster
\item"(i)" $\forall\sigma>0\ \mu\big(\{y\in M:
\ |y(\tau)-\tau|\geqslant\sigma\}\big)\to 0$, $\tau\to\pm\infty$;
\item"(ii)" $\mu\big(\{y\in M:\ y(\tau)-\tau,\tau-y(\tau)\text{ are bounded
on }\RR\}\big)=0$.
\endroster
\endproclaim

\remark{Remark}
The gamma noise measure is a concrete example of generalized white noise
measure. The latter measures were introduced in [5], see also [1,4]. The
first statement of Theorem 1 is valid for all the generalized white noise
measures (it is a simple consequence of L\'evy--Khinchine representation).
The second one holds for compound Poisson measures such that the
measure $\beta$ in L\'evy triple satisfies the following conditions: $\beta$
is an infinite measure on $(0,\infty)$; $\forall\delta>0$
$\beta\big((\delta,\infty)\big)<\infty$; $m_1(\beta)=\int_0^\infty
ud\beta(u)<\infty$. The third statement holds for compound Poisson
measures, satisfying additional requirement:  $m_1(\beta)=1$.

\endremark

2. It follows from \thetag1 that the gamma noise measure $\mu$ is analytic
one on $(N',\Cal B(N'))$. Then the function
$$
e_\mu(x;\theta)=\exp\{\l x,\theta\r+\l 1,\log(1-\theta)\r\},
\ x\in N',\ \theta\in N
$$
generates the system $\bold P^\mu$ of Appell polynomials, corresponding to
$\mu$, i.e. $e_\mu(x;\theta)=\sum_{n=0}^\infty
\frac{\l P^\mu_n(x),\theta^{\otimes n}\r}{n!}$ (see [12]). Let a
vector-function $\NC\ni\theta\mapsto\alpha(\theta)\in\NC$ be holomorphic
at zero. Then
$$
e_\mu^\alpha(x;\theta):=e_\mu(x;\alpha(\theta)),\ x\in N',\ \theta\in N
$$
is generating function of the system $\bold P^{\mu,\alpha}$ of generalized
Appell polynomials, corresponding to measure $\mu$ and function $\alpha$
(for more details see [11], in one-dimensional case $N'=\H=N=\RR$ see [3]).

Let $\pi$ be Poisson measure on $(N',\Cal B(N'))$. Recall that the
generalized Appell polynomials $\bold P^{\pi,\alpha_\pi}$ corresponding to
$\pi$ and $\alpha_\pi(\theta)=\log(1+\theta)$ coincide with Charlier
polynomials, that are orthogonal with respect to $\pi$, see [14,15]. In the
case of gamma noise measure the following statement is true (see [7,8]).

\proclaim{Theorem 2}
Let $\alpha_\mu(\theta)=\frac\theta{\theta-1}$, $\theta\in N$. The
generalized Appell polynomials $\bold P^{\mu,\alpha_\mu}$ are orthogonal
with respect to gamma noise measure $\mu$.
\endproclaim

In one-dimensional case $\bold P^{\mu,\alpha_\mu}$ coincides with the
system of Laguerre polynomials (see, e.g., [3]). So we can call the
polynomials from the system $\bold J^\mu:=\bold P^{\mu,\alpha_\mu}$
infinite-dimensional Laguerre polynomials. As in Gaussian and Poissonian
analysis one can construct the rigging $(N)^{-1}_{\bold J^\mu}\supseteq
\roman L_2(N',\mu)\supseteq(N)^1_{\bold J^\mu}$, using the system
$\bold J^\mu=\{\l J^\mu_n(\cdot),\varphi^{(n)}\r,
\ \varphi^{(n)}\in\NC^{\widehat\otimes n},\ n\in\Bbb Z_+\}$
($\widehat\otimes$ denotes symmetric tensor product)
of Laguerre polynomials. Namely, $(N)^1_{\bold J^\mu}=\prlim_{p,q\in\N}
(\H_p)^1_q$, $(N)^{-1}_{\bold J^\mu}=\indlim_{p,q\in\N}(\H_{-p})^{-1}_{-q}$,
where $(\H_p)^1_q$ is a completion of the space of continuous polynomials
on $N'$
$\Cal P(N')=\big\{\varphi(x)=\sum_{k=0}^n\l J^\mu_k(x),\varphi^{(k)}\r$,
$\varphi^{(k)}\in\NC^{\widehat\otimes k}$,
$k=\overline{0,n},\ n\in\Bbb Z_+\big\}$ with respect to the Hilbert norm
$\|\varphi\|_{p,q,1}^2=\sum_{n=0}^\infty (n!)^22^{qn}|\varphi^{(n)}|_p^2$,
$\varphi\in\Cal P(N')$, $(\H_{-p})^{-1}_{-q}$ is dual of $(\H_p)^1_q$ with
respect to $\roman L_2(N',\mu)$.

On another hand, using $\bold J^\mu$ one can prove a nondegeneracy of the
gamma noise measure $\mu$ and construct the dual Appell system $\bold
Q^\mu$ of distributions on $N'$ and the rigging $(N)^{-1}_\mu\supset\roman
L_2(N',\mu)\supset(N)^1$ of $\roman L_2(N',\mu)$ by Kondratiev spaces of
distributions and test functions respectively (see [12]). As in the case of
Poisson measure $\pi$ the equalities $(N)^{-1}_\mu=(N)^{-1}_{\bold J^\mu}$,
$(N)^1=(N)^1_{\bold J^\mu}$ are valid, so we have a description of
$(N)^{-1}_{\bold J^\mu}$ and $(N)^1_{\bold J^\mu}$ via characterization
theorems from [12].

We apply this result for solving of the gamma counterpart of Verhulst type
equation
$$
Y_t=Y_0+r\int_0^tY_s\diamond(1-Y_s)\roman ds
+a\int_0^tY_s\diamond(1-Y_s)\diamond\xi_s'\roman ds, \tag3
$$
modelling the population growth in a crowded stochastic environment.
Here $\{Y_t$, $t\in\RR_+\}\subset(N)^{-1}_\mu$ is unknown generalized
stochastic process, $r\in\RR$, $a>0$, $\diamond$ means Wick product (for
more details on this and more general Wick--Skorokhod equations see [9],
Poisson counterpart of \thetag3 was considered in [15]).

\proclaim{Theorem 3}
Let $Y_0\in(N)^{-1}_\mu$ be such that $\L Y_0,1\R\ne 0$
($\L\cdot,\cdot\R$ denotes the dual pairing between $(N)^{-1}_\mu$ and
$(N)^1$, given by inner product in $\roman L_2(N',\mu)$). Then
$$
Y_t(x)=\big[1+\big(Y_0^{\diamond(-1)}-1\big)\diamond
e_\mu^{\alpha_\mu}(-r-ax;\chi_{[0,t]})\big]^{\diamond(-1)},\ x\in N'
$$
is the unique $(N)^{-1}_\mu$-solution to \thetag3.
\endproclaim

The proof uses $S_\mu^{\alpha_\mu}$-transform of $\Phi\in(N)^{-1}_\mu$:
$(S_\mu^{\alpha_\mu}\Phi)(\theta)=\L\Phi,e_\mu^{\alpha_\mu}(\cdot;\theta)\R$,
being the analogue of $S$-transform of Gaussian and $S_\pi$-transform of
Poissonian analysis. Namely, we have for any $\Phi,\Psi\in(N)^{-1}_\mu
=(N)^{-1}_{\bold J^\mu}$ that $S_\mu^{\alpha_\mu}(\Phi\diamond\Psi)(\theta)
=S_\mu^{\alpha_\mu}(\Phi)(\theta)S_\mu^{\alpha_\mu}(\Psi)(\theta)$,
$\theta\in N$. So $S_\mu^{\alpha_\mu}$ transfers the equation \thetag3 to
deterministic one as well as $S$ and $S_\pi$. It remains to apply the
standard arguments of white noise functional approach to Wick--Skorokhod
equations [9].

\bigskip
{\smc
A. V. Gorbunov, Assistant, Department of Economics, Kiev University,
Volody\-myr\-ska st., 64, Kiev, 252033, Ukraine.}
\bigskip
{\smc G. F. Us, Ph.D., Associate Professor, Department of Mechanics and
Mathematics, Kiev University, Volodymyrska st., 64, Kiev, 252033, Ukraine.}

\enddocument